\documentclass[journal]{IEEEtran}
\usepackage{amsmath,amssymb,amsfonts}
\usepackage{graphicx}
\usepackage{xcolor}
\usepackage{enumerate}
\usepackage{amssymb}
\usepackage{lipsum}

\newtheorem{proposition}{Proposition}[section]
\newtheorem{lemma}{Lemma}[section]

\newtheorem{remark}{Remark}[section]

\newtheorem{assumption}{Assumption}[section]
\makeatletter
\DeclareRobustCommand{\qed}{%
  \ifmmode 
  \else \leavevmode\unskip\penalty9999 \hbox{}\nobreak\hfill
  \fi
  \quad\hbox{\qedsymbol}}
\newcommand{\openbox}{\leavevmode
  \hbox to.77778em{%
  \hfil\vrule
  \vbox to.675em{\hrule width.6em\vfil\hrule}%
  \vrule\hfil}}
\newcommand{\qedsymbol}{\openbox}
\newenvironment{proof}[1][\proofname]{\par
  \normalfont
  \topsep6\p@\@plus6\p@ \trivlist
  \item[\hskip\labelsep\itshape
    #1.]\ignorespaces
}{%
  \qed\endtrivlist
}
\newcommand{\proofname}{Proof}
\makeatother

\begin{document}

\title{Backstepping Control of a Hyperbolic PDE System with Zero Characteristic Speed States}

\author{Gustavo A. de Andrade, Rafael Vazquez, Iasson Karafyllis, and Miroslav Krstic
\thanks{Gustavo A. de Andrade, is with the Departament of Automation and Systems, Universidade Federal de Santa Catarina, 88040-370 Florian\'{o}polis, Brazil (e-mail: gustavo.artur@ufsc.br). }
\thanks{Rafael Vazquez is with the Department of Aerospace Engineering, Universidad de Sevilla, 41092 Sevilla, Spain (e-mail: rvazquez1@us.es).}
\thanks{Iasson Karafyllis is with the Department of Mathematics, National Technical University of Athens, 15780 Athens, Greece (iasonkar@central.ntua.gr).}
\thanks{Miroslav Krstic is with the Department of Mechanical and Aerospace
Engineering, University of California San Diego, La Jolla, CA 92093-
0411 USA (e-mail: krstic@ucsd.edu).}}
\markboth{}%
{Shell \MakeLowercase{\textit{et al.}}: Bare Demo of IEEEtran.cls for IEEE Journals}

\maketitle

\begin{abstract}
While for coupled hyperbolic PDEs of first order there now exist numerous PDE backstepping designs, systems with zero speed, i.e., without convection but involving infinite-dimensional ODEs, which arise in many applications, from environmental engineering to lasers to manufacturing, have received virtually no attention. In this paper, we introduce single-input boundary feedback designs for a linear 1-D hyperbolic system with two counterconvecting PDEs and $n$ equations (infinite-dimensional ODEs) with zero characteristic speed. The inclusion of zero-speed states, which we refer to as {\em atachic}, may result in non-stabilizability of the plant. We give a verifiable condition for the model to be stabilizable and design a full-state backstepping controller which exponentially stabilizes the origin in the $\mathcal{L}^{2}$ sense. In particular, to employ the backstepping method in the presence of atachic states, we use an invertible Volterra transformation only for the PDEs with nonzero speeds, leaving the zero-speed equations unaltered in the target system input-to-state stable with respect to the decoupled and stable counterconvecting nonzero-speed equations. Simulation results are presented to illustrate the effectiveness of the proposed control design.
\end{abstract}

\IEEEpeerreviewmaketitle

\section{Introduction}

In the last two decades,  boundary stabilization of hyperbolic partial differential equations (PDEs) has been widely studied in the literature due to its applicability in several processes, such as open water channels \cite{deHalleux2003}, traffic flow \cite{Espitia2022}, flow through pipelines \cite{Gugat2011}, oil wells \cite{Aamo2013} and electrical transmission lines \cite{Bastin2016}. The area has achieved a  stage that is rather advanced (``mature'' is the adjective that might come to the creativity-challenged mind) thanks to the power of the PDE backstepping method in the design of control laws and observers. Results covering the  stabilization of $(n+m)\times(n+m)$  hyperbolic  systems consisting of $n$ (uncontrolled) equations convecting in one direction and $m$ controlled equations (counter-)convecting in the opposite direction were introduced in \cite{hu2016}. The extension of this approach for the output-feedback regulation with additional disturbances is proposed in \cite{deutscher2018}. An adaptive observer design for $(n+1)\times(n+1)$ hyperbolic systems can be found in \cite{anfinsen2016}, where the methodology can be used even in the case of unknown or incorrect parameters. Recently, the control design of hyperbolic systems with nonstrict-feedback couplings with ordinary differential equations (ODEs) has also been studied, see e.g. \cite{deAndrade2018} and references therein. All these developments assume nonzero characteristic speeds, becoming inapplicable otherwise.

Outside of the backstepping literature, there exists a few studies on exact boundary controllability for very specific hyperbolic systems with identically zero or vanishing characteristic speeds, as can be seen in \cite{Karafyllis2018,coron2009} or \cite{Tatsien2009}. An approach based on static output feedback controllers was proposed in \cite{yong2019boundary} but requires the so-called structural stability conditions which are rather strict on the coefficients of the system.

We refer to hyperbolic systems containing states with zero velocity as {\em atachic} (meaning `with no velocity' in Greek; recall the term {\em isotachic} in \cite{hu2016} and the special properties of PDEs of equal characteristic speeds). The disregard in the control literature for such systems does not imply they are not of practical interest. On the contrary, multiple applications do exist. A first example is a model of heat transfer dynamics in solar thermal plants based on  direct steam generation technology \cite{guo2017}. In this application, solar radiation heats the water to generate superheated steam which is used by a turbine generator to convert the thermal energy into electricity. The pipe temperature dynamics is usually described by a hyperbolic PDE with zero characteristic speed, which is coupled with the multiphase flow equations. The exponential stabilization of the states and reference tracking is crucial to correct and safe operation of the system, and can be performed by manipulating the boundary mass flow rate, using a pump placed at the inlet of the solar field. Addiabatic flows, such as the Saint-Venant equations or the isentropic and full Euler equation for gas dynamics  \cite{coron2009}, which are of practical interest on accounting of the trend to perate combined sewer systems and other channel networks,  also admit zero characteristic speed.

Other application that fits into the atachic framework is the intensity dynamics of the laser beam \cite{ren2012laser}. In several laser applications, the maximization of the energy extracted from the laser pulse is critical for the process eficience. For example, in the polysilicon process for manufacturing flat panel displays, one of the main problems is to obtain enough instantaneous laser power to melt as large an area as desired. Photolithography is another example where the optical exposures must be accomplished with fewer pulses of higher energy.

A few more examples include models with thermoacoustic instabilities --- a zero transport velocity in thermoacoustics is a direct consequence of the second law of thermodynamics --- double-pass laser amplifiers~\cite{ren2012laser}, neurofilament transport in axons~\cite{craciun2005dynamical}, and biomass production in photobioreactors \cite{fernandez2016}.

Motivated by these applications, this paper aims to extend  the infinite-dimensional backstepping methodology to what we denote, extending the $(n+m)\times(n+m)$ notation, as  $(1+n+1)\times(1+n+1)$ 1-D hyperbolic systems, which contain one rightward convecting unactuated state, $n$ non-convecting/zero-speed/atachic unactuated states, and one leftward-convecting state with boundary actuation. 

As explained, previous results such as \cite{hu2016} or \cite{dimeglio2013} are inapplicable since they would result in a controller with infinite gain. Additionally, it is shown that not all $(1+n+1)\times(1+n+1)$ systems are stabilizable. The homogeneous part of the zero-speed equation must be asymptotically stable for the overall system to be stabilizable without other restrictions. 

We apply a backstepping transformation only to the PDEs with nonzero speeds, leaving the state of the zero-speed equation unaltered in the target system, but making the target zero-speed equation input-to-state stable with respect to the decoupled and stable counterconvecting nonzero-speed target PDEs.  Compared with other results in the literature for hyperbolic PDEs containing states with zero characteristic speeds, such as \cite{Tatsien2009}, our approach can be applied to a  richer family of hypebolic systems that can be unstable in the nonzero speed part of the plant. In particular we provide numerical simulations to show the effectiveness of the method for an open-loop unstable case. Part of these contributions was previously published in preliminary form in the conference paper \cite{deAndrade2022} for a $(1+1+1)\times(1+1+1)$ system.

The rest of the paper is organized as follows. In Section \ref{section:problem_statement}, we present the control problem and some properties of the equations that motivate our assumptions on controlling the system. In Section \ref{section:control_design}, we design a stabilizing control law using the backstepping methodology. The results are illustrated using numerical simulations in Section \ref{section:simulations}. Finally, Section \ref{section:conclusions} provides some concluding remarks and directions of future work.

\subsection*{Notation}

For a given $u:\mathbb{R}_{+}\times[0,1]\rightarrow\mathbb{R}$, we use the notation $u[t]$ to denote the profile at certain $t\geq 0$, i.e., $(u[t])(x)=u(t,x)$ $\forall x\in[0,1]$. $\mathcal{L}^{2}(0,1)$ denotes the set of equivalence classes of measurable functions $f:[0,1]\rightarrow \mathbb{R}$ for which $\|f\|_{2}=\left(\int_{0}^{1}|f(x)|^{2}dx \right)^{1/2}<+\infty$. For an interval $I\subset\mathbb{R}_{+}$, the space $\mathcal{C}^{0}(I;\mathcal{L}^{2}(0,1))$ is the space of continuous mappings $I\ni t \rightarrow u[t]\in\mathcal{L}^{2}(0,1)$. Finally, $\mathcal{H}^{1}(0,1)$ denotes the Sobolev space of functions in $\mathcal{L}^{2}(0,1)$ with all its first-order weak derivatives in $\mathcal{L}^{2}(0,1)$.

\section{Problem Statement}\label{section:problem_statement}

For $n\geq 1$, consider the following set of $(1+n+1)$ hyperbolic system:
\begin{align}
    \partial_t u(t,x)=&- \lambda_1\partial_x u (t,x)+\sigma_{12}p(t,x) + \Theta_{1}v(t,x), \label{eq_u}\\
        \partial_t v(t,x)=&\Omega_{1}u(t,x)+ \Omega_{2}p(t,x) + \Psi v(t,x),\label{eq_vn}\\
    \partial_t p(t,x) =& \lambda_2\partial_{x} p(t,x)+\sigma_{21}u(t,x) + \Theta_{2}v(t,x),\label{eq_p}\\
    u(t,0) =& U(t) + qp(t,0), \label{bc_u}\\
     p(t,1) =& \rho u(t,1),\label{bc_p}
\end{align}
\noindent where $t\in[0,\infty)$ is the time, $x\in[0,1]$ is the space, the states are given by $u$, $p$ and $v = \left(v_{1}, \;\;\; \dots \;\;\; , v_{n} \right)$, and the control action is $U$. The transport speeds satisfy $\lambda_{1}>0>-\lambda_{2}$, and $\rho,q$ are nonzero reflection coefficients. The other coefficients of the system are 
\begin{align*}
    \Theta_{1} &= \left(\theta_{11} \;\;\; \dots \;\;\; \theta_{1n} \right),
    & \Theta_{2} &= \left(\theta_{21} \;\;\; \dots \;\;\; \theta_{2n} \right),\\
    \Omega_{1} &= \left(\omega_{11} \;\;\; \dots \;\;\; \omega_{n1} \right)^{T}, 
    & \Omega_{2} &= \left(\omega_{12} \;\;\; \dots \;\;\; \omega_{n2} \right)^{T},\\
    \Psi &=\{\psi_{ij}\}_{1\leq i\leq n, 1\leq j\leq n}. 
\end{align*}

Finally, the initial condition of \eqref{eq_u}-\eqref{bc_u} are
\begin{align}
    u(0,x) = u_{0}(x),\; v(0,x)=v_{0}(0,x),\; p(0,x)=p_{0}(x),\label{eq_ini}
\end{align}
\noindent with $u_{0}\in\mathcal{L}^{2}(0,1)$, $v_{0}\in\left(\mathcal{L}^{2}(0,1)\right)^{n}$ and $p_{0}\in\mathcal{L}^{2}(0,1)$.

\begin{figure}[b!]
\centering
\includegraphics[width=60mm]{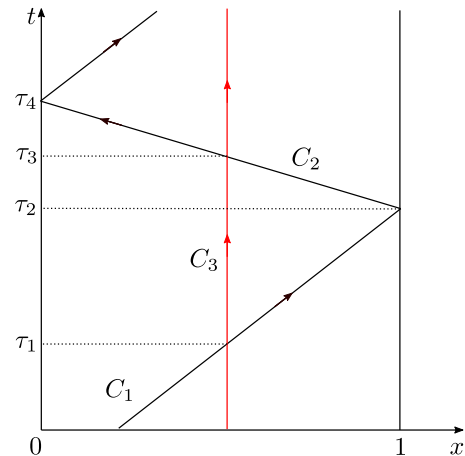}
\caption{Characteristic lines of system. The characteristic lines $C_{1}$ (with slope $\lambda_{1}$) and $C_{2}$ (with slope $-\lambda_{2}$) correspond to \eqref{eq_u} and \eqref{eq_p}, respectively, whereas $C_{3}$ corresponds to \eqref{eq_vn} with $n=1$. The reflection mechanism is illustrated at the points $x=0$ and $x=1$ at the time instants $\tau_{2}$ and $\tau_{4}$, respectively.}\label{fig:system_characteristics}
\end{figure}

System \eqref{eq_u}-\eqref{eq_ini} is hyperbolic with two characteristic speeds with opposite signs, associated with \eqref{eq_u}-\eqref{eq_p}, respectively, and $n$ identically zero speeds for \eqref{eq_vn}. The latter means that the characteristics corresponding to \eqref{eq_vn} are vertical on the $(x,\,t)$ plane (see Figure \ref{fig:system_characteristics}). As shown in \cite{Tatsien2009}, the stabilizability of \eqref{eq_u}-\eqref{eq_ini} can be obtained (without any constraints on the parameters) by imposing boundary controllers and an in-domain controller for the equation with zero characteristic speed. However if internal controllers are not realizable, then the stabilizability can be obtained only for very strict cases, as will be shown in the next section. For that reason, we state the following assumption:
\begin{assumption}\label{assumption:psi_matrix}
It is assumed that $\Psi$ is a Hurwitz matrix.
\end{assumption}

\subsection{Stabilizability of systems with zero characteristic speeds}

In this section we will show, with a simplified version of system \eqref{eq_u}-\eqref{eq_ini}, that zero characteristic speeds may result in non stabilizability of the system. Although this example  does not consider counter convective PDEs, it is to be expected that this result is valid for the general form \eqref{eq_u}-\eqref{eq_ini}, but accompanied by more complex expressions.  Before the formal result is presented, we will discuss an example to build insights of the problem. 

Consider the following hyperbolic system:
\begin{align}
    \partial_t u(t,x) =& -\partial_{x} u(t,x),\label{eq:simplified1_u}\\
    \partial_t v(t,x) =& \psi v(t,x) + u(t,x), \label{eq:simplified1_v}\\
    u(t,0) =& U(t).\label{eq:bc_simplified}
\end{align}
\noindent with $\psi>0$, and $u(0,x)=v(0,x)=0$ for all $x\in[0,1]$.

One way to argue the lack of stabilizability of \eqref{eq:simplified1_u}-\eqref{eq:bc_simplified} is to convert the equations into a form in which a part of the dynamics is autonomous and unstable. To do that, first note that after $t\geq 1$ the solution of \eqref{eq:simplified1_u} is $u(t,x) = U(t-x)$, and thus, the atachic equation \eqref{eq:simplified1_v} can be rewritten to
\begin{align}
    \partial_t v(t,x)= \psi v(t,x) + U(t-x).    
\end{align}

Now, set any $0<x_{1}\leq 1$, and define $v_{1}(t) = v(t,0)$ and $v_{2}(t)=v(t, x_{1})$. Then,
\begin{align}
    \dot{v}_{1}(t) &= \psi v_{1}(t) + U(t),\label{eq:ode_v1}\\
    \dot{v}_{2}(t) &= \psi v_{2}(t) + U(t-x_{1}).\label{eq:ode_v2}
\end{align}

Importantly, note that the history of $U(t)$ for $t\in[-1,0]$ for this delay system are needed since $0< x_{1}\leq 1$ (these values will come from the initial condition of $v$ in (\ref{eq:simplified1_u})).

For $t\geq x_{1}$, define $w(t)=v_{1}(t-x_{1})-v_{2}(t)$. Then, it follows that
\begin{align}
    \dot{w}(t) = \psi w(t).\label{eq:ode_w}
\end{align}

Note that \eqref{eq:ode_w} is autonomous and unstable, and thus, unless $w(x_{1})=0$, we have $|w(t)|\rightarrow \infty$ as $t\rightarrow\infty$.  By definition, the only way to have $w(x_{1}) = 0$ is if $v_{2}(x_{1})=v_{1}(0)$. Thus using the explict solution of \eqref{eq:ode_v2}, replacing the expressions of $v_{1}$ and $v_{2}$, and writting $x$ instead of $x_{1}$, we have that
\begin{align}
   v_{0}(0) = \mathrm{e}^{\psi x}v_{0}(x)+\int_{-x}^{0}\mathrm{e}^{-\psi \tau}U(\tau)d\tau
\end{align}
\noindent must be verified by $U(t)$ for $t\in[-1,0]$ so that the system can be stabilized. If that is the case, one has $w(t)=0$ for all $0<x_{1}\leq 1$ and therefore it holds that for $t\geq x$, $v(t,x)=v(t-x,0)=v_{1}(t-x)$, where $v_{1}$ satisfies \eqref{eq:ode_v1}, and so designing a stabilizing control law for $v_{1}$ also stabilizes for $v$.

This quick argument shows that one cannot expect stabilizability of (\ref{eq:simplified1_u})--(\ref{eq:bc_simplified}) if $\psi >0$ (it would also hold for $\psi =0$).

We now formalize and generalize the above arguments and prove the lack of stabilizability of the following  hyperbolic system:
\begin{align}
    \partial_t u(t,x) =& -\lambda \partial_{x} u(t,x),\label{eq:gen_simplified_u}\\
    \partial_t v(t,x)=& \psi v(t,x) + \omega u(t,x), \label{eq:gen_simplified_v}\\
    u(t,0) =& U(t),\label{eq:gen_bc_simplified}
\end{align}
\noindent where $\lambda >0$, $\psi \geq0$ and $\omega \in \mathbb{R}$.  

The initial condition of \eqref{eq:gen_simplified_u}-\eqref{eq:gen_bc_simplified} is
\begin{align}
    u(0,x) &= u_{0}(x), & v(0,x) & = v_{0}(x)\label{eq:gen_ini_simplified}
\end{align}
\noindent with $u_{0}\in\mathcal{H}^{1}(0,1)$ and $v_{0}\in\mathcal{H}^{1}(0,1)$.

The idea consists in proving the existence of a continuous functional $P(u,v)$, with  $P(u,v)=0$ for which the set $S\{(u,v):W(u,v)>0\}$ is non-empty and $\frac{d W}{dt}(u,v)\geq 0$ for all $(u,v)\in S$ and for all $U(t)\in \mathbb{R}$.

\begin{proposition}\label{prop:1}
    Let $\lambda >0$, $\psi \geq 0$, and $\omega \in \mathbb{R}$ be given constants. Let $S\subset \mathcal{H}^{1}(0,1)\times \mathcal{H}^{1}(0,1)$ be the linear subspace
    \begin{align}
        S = \{(u,v)\in \mathcal{H}^{1}(0,1)\times \mathcal{H}^{1}(0,1): P(u,v)=0 \}\label{eq:linear_subspace}
    \end{align}
    \noindent where $P:\mathcal{H}^{1}(0,1)\times \mathcal{H}^{1}(0,1)\rightarrow \mathcal{H}^{1}(0,1)$ is the linear operator defined by
    \begin{multline}
        (P(u,v))(x) = v(x) - \mathrm{e}^{-\lambda ^{-1}\psi x}v(0) \\
        +\lambda ^{-1}\omega \int_{0}^{x}\exp^{-\lambda^{-1}\psi(x-s)}u(s)ds,\label{eq:operator_P}
    \end{multline}
    \noindent for all $(u,v)\in \mathcal{H}^{1}(0,1)\times \mathcal{H}^{1}(0,1)$ and $x\in[0,1]$.

    Then the following property hold for all $(u_{0},v_{0})\in \mathcal{H}^{1}(0,1)\times \mathcal{H}^{1}(0,1)\setminus S$:

    \noindent\textbf{(P)} For every input $U\in \mathcal{C}^{1}(\mathbb{R}_{+})$ with $U(0)=u_{0}(0)$ the corresponding unique solution $(u,v)\in \left(\mathcal{C}^{0}(\mathbb{R}_{+};\mathcal{H}^{1}(0,1))\cap\right.$ $\left. \mathcal{C}^{1}(\mathbb{R}_{+};\mathcal{L}^{2}(0,1))\right)^{2}$ of the initial-boundary value problem \eqref{eq:gen_simplified_u}-\eqref{eq:gen_ini_simplified} does not satisfy $\lim_{t\rightarrow\infty}(\|u[t]\|_{2}) = \lim_{t\rightarrow\infty}(\|v[t]\|_{2}) = 0.$ 
\end{proposition}

\begin{proposition}\label{prop:2}
    Let $\lambda >0$, $\psi \geq 0$, $\omega \in \mathbb{R}$ be given constants. Then the linear subspace $S\subset \mathcal{H}^{1}(0,1)\times \mathcal{H}^{1}(0,1)$ defined by \eqref{eq:linear_subspace} is positively invariant. In other words, for every $(u_{0},v_{0})\in S$ and for every input $U\in \mathcal{C}^{1}(\mathbb{R}_{+})$, with $U(0)=u_{0}(0)$, the corresponding unique solution $(u,v)\in \left(\mathcal{C}^{0}(\mathbb{R}_{+};\mathcal{H}^{1}(0,1))\cap\right.$ $\left. \mathcal{C}^{1}(\mathbb{R}_{+};\mathcal{L}^{2}(0,1))\right)^{2}$ of the initial value problem \eqref{eq:gen_simplified_u}-\eqref{eq:gen_ini_simplified} satisfies $(u[t],v[t])\in S$ for all $t\geq 0$.
\end{proposition}

\begin{proposition}\label{prop:3}
    Let $\lambda >0$, $\psi \geq 0$, $\omega \in \mathbb{R}$ be given constants. Let $S\subset \mathcal{H}^{1}(0,1)\times \mathcal{H}^{1}(0,1)$ be the linear subspace defined by \eqref{eq:linear_subspace}. Then the following property holds for all $(u_{0},v_{0})\in \mathcal{H}^{1}(0,1)\times \mathcal{H}^{1}(0,1)\setminus S$:

    \noindent \textbf{(P')} For every input $U\in \mathcal{C}^{1}(\mathbb{R}_{+})$ with $U(0)=u_{0}(0)$ the corresponding unique solution $(u,v)\in \left(\mathcal{C}^{0}(\mathbb{R}_{+};\mathcal{H}^{1}(0,1))\cap\right.$ $\left. \mathcal{C}^{1}(\mathbb{R}_{+};\mathcal{L}^{2}(0,1))\right)^{2}$ of the initial-boundary value problem \eqref{eq:gen_simplified_u}-\eqref{eq:gen_ini_simplified} does not satisfy $\lim_{t\rightarrow\infty}v(t,x) = 0.$ for all $x\in[0,1]$. 
\end{proposition}

The proofs of Propositions~\ref{prop:1},~\ref{prop:2} and \ref{prop:3} are given in the Appendix.
\begin{remark}
    Property \textbf{(P')} is different from property \textbf{(P)} since  \textbf{(P')} deals with pointwise convergence to zero while property  \textbf{(P)} deals with convergence to zero in $\mathcal{L}^{2}$. It is well-known that convergence in $\mathcal{L}^{2}$ does not imply pointwise convergence and pointwise convergence dos not imply convergence in $\mathcal{L}^{2}$. Moreover, property \textbf{(P')} deals with the pointwise convergence of $v$ only, while \textbf{(P)} deals with convergence in $\mathcal{L}^{2}$ of $(u,v)$.
\end{remark}

\begin{remark}
    System \eqref{eq:gen_simplified_u}-\eqref{eq:gen_bc_simplified} is equivalent to the system
    \begin{align}
        \partial_{t} u(t,x) &= -\lambda\partial_{x}u(t,x),\\
        \partial_{t}w(t,x) &= \psi w(t,x),\\
        \dot{y}(t) &= \psi y(t) + \omega U(t),\label{eq:ode_y}\\
        u(0) &= U(t),\\
        w(t,0) &= 0,
    \end{align}
    \noindent with state $(u,w,y)\in \mathcal{H}^{1}(0,1)\times \mathcal{H}^{1}(0,1)\times \mathbb{R}$. To see this, notice that transformation \eqref{eq:definition_w} in Appendix \ref{app:proof_prop2} and the definition $y(t) = v(t,0)$ give us \eqref{eq:ode_y}. On the other hand, the inverse transformation is given by $v(x) = w(x)+\mathrm{e}^{-\lambda^{-1}\psi x}y - \lambda^{-1}\omega\int_{0}^{x}\mathrm{e}^{-\lambda^{-1}\psi(x-s)}u(s)ds$ for $x\in[0,1]$ and transforms \eqref{eq:ode_y} to \eqref{eq:gen_simplified_u}-\eqref{eq:gen_bc_simplified}. The system may be defined on the invariant subspace $S$ defined by \eqref{eq:linear_subspace} and corresponds to the case $w=0$ for \eqref{eq:ode_y}. Thus system \eqref{eq:gen_simplified_u}-\eqref{eq:gen_bc_simplified} with $\omega\neq 0$ on the invariant subspace $S$ may be stabilized by the feedback law
    \begin{align*}
        U(t) = -\frac{k}{\omega}v(0),
    \end{align*}
    where $k>\omega$ is a design constant. When $\omega=0$,  \eqref{eq:gen_simplified_u}-\eqref{eq:gen_ini_simplified} can be stabilized on the subspace $S' = \{(u,v)\in S: v(0)=0 \}$ with $U(t)=0$.
\end{remark}

\section{Control Design}\label{section:control_design}

Having justified Assumption~\ref{assumption:psi_matrix}, in this section we design a backstepping controller so that the null solution of \eqref{eq_u}-\eqref{eq_ini} becomes stable. It will be assumed that the full-state measurements are available for the control law.

Consider the following Volterra transformation
    \begin{align}
    \alpha =& u - \int_{x}^{1}K_{1}(x,\xi)u(t,\xi)d\xi  -\int_{x}^{1}K_{2}(x,\xi)p(t,\xi)d\xi \nonumber\\
    &-\int_{x}^{1}G(x,\xi)v(t,\xi)d\xi,\label{eq_transf1}\\
    \beta =& p - \int_{x}^{1}Q_{1}(x,\xi)u(t,\xi)d\xi - \int_{x}^{1}Q_{2}(x,\xi)p(t,\xi)d\xi \nonumber\\
    &-\int_{x}^{1}R(x,\xi)v(t,\xi)d\xi,\label{eq_transf2}
\end{align}
\noindent where the kernels $K_{i}$, $Q_{i}$, for $i\in \{1,2\}$, $G = (G_{1}\;\;\;\hdots \;\;\; G_{n})$ and $R=(R_{1}\;\;\;\hdots \;\;\; R_{n})$ satisfy the following PDEs:
\begin{align}
    \lambda_{1}\partial_{x} K_{1} + \lambda_{1}\partial_{\xi}K_{1} =& -\sigma_{21}K_{2} - G\Omega_{1},\label{eq:pde1_kernels}\\
    \lambda_{1}\partial_{x} K_{2} - \lambda_{2}\partial_{\xi}K_{2} =&-\sigma_{12}K_{1} -  G\Omega_{2},\label{eq:pde2_kernels}\\
    \lambda_{1}\partial_{x} G =&  -K_{1}\Theta_{1} - K_{2}\Theta_{2} -  G\Psi,\label{eq:pde3_kernels}\\
    \lambda_{2}\partial_{x} Q_{1} - \lambda_{1}\partial_{\xi}Q_{1} =&\sigma_{21}Q_{2} + R\Omega_{1},\label{eq:pde4_kernels}\\
    \lambda_{2}\partial_{x} Q_{2} + \lambda_{2}\partial_{\xi}Q_{2} =& \sigma_{12}Q_{1} + R\Omega_{2},\label{eq:pde5_kernels}\\
    \lambda_{2}\partial_{x} R =& Q_{1}\Theta_{1} + Q_{2}\Theta_{2} + R\Psi,\quad\label{eq:pde6_kernels}
\end{align}
\noindent with boundary conditions
\begin{align}
    K_{1}(x,1) &= \frac{\lambda_{2}\rho}{\lambda_{1}}K_{2}(x,1), & Q_{1}(x,x) &= \frac{\sigma_{21}}{\lambda_{1}+\lambda_{2}},\label{eq:bc1_kernels}\\
    K_{2}(x,x) &= -\frac{\sigma_{12}}{\lambda_{1}+\lambda_{2}}, & Q_{2}(x,0) &= \frac{\lambda_{1}}{\lambda_{2}q} Q_{1}(x,0),\label{eq:bc2_kernels}\\
    G(x,x) &=  -\frac{\Theta_{1}}{\lambda_{1}}, & R(x,x) &= \frac{\Theta_{2}}{\lambda_{2}}.\label{eq:bc3_kernels}
\end{align}
\noindent for $k\in\{1,\hdots , n\}$.

These kernels evolve in the triangular domain 
\begin{align*}
    \mathcal{T} = \{(x,\xi)\in \mathbb{R}^{2}:\; 0\leq x \leq \xi \leq 1\}.
\end{align*}

The explicit solution of \eqref{eq:pde3_kernels} and \eqref{eq:pde6_kernels}, together with the boundary condition \eqref{eq:bc3_kernels}, is given by
\begin{align}
    G(x,\xi) =& -\frac{1}{\lambda_{1}}\Theta_{1}\Phi_{1}(x,\xi),\nonumber\\
    &-\frac{1}{\lambda_{1}}\int_{\xi}^{x}K_{1}(\tau,\xi)\Theta_{1}\Phi_{1}(x,\tau)d\tau \nonumber\\ 
    &-\frac{1}{\lambda_{1}}\int_{\xi}^{x}K_{2}(\tau,\xi)\Theta_{2}\Phi_{1}(x,\tau)d\tau,\label{eq:solutionG}\\
    R(x,\xi) =& \frac{1}{\lambda_{2}}\Theta_{2}\Phi_{2}(x,\xi) \nonumber\\ 
    &+\frac{1}{\lambda_{2}}\int_{\xi}^{x}Q_{1}(\tau,\xi)\Theta_{1}\Phi_{2}(\tau,x)d\tau \nonumber\\ 
    &+\frac{1}{\lambda_{2}}\int_{\xi}^{x}Q_{2}(\tau,\xi)\Theta_{2}\Phi_{2}(\tau,x)d\tau.\label{eq:solutionR}
\end{align}
\noindent where $\Phi_{k}(x,\xi) = \mathrm{e}^{\Psi\frac{(\xi-x)}{\lambda_{k}}}$, for $k\in\{1,2\}$, is the state-transition matrix.

Plugging \eqref{eq:solutionG} and \eqref{eq:solutionR} into \eqref{eq:pde1_kernels}-\eqref{eq:pde2_kernels} and \eqref{eq:pde4_kernels}-\eqref{eq:pde5_kernels}, respectively, one obtain the following system of integro-differential equations: 
\begin{align}
    &\lambda_{1}\partial_{x} K_{1} + \lambda_{1}\partial_{\xi}K_{1} = -\sigma_{21}K_{2} + \lambda_{1}^{-1}\Theta_{1}\Phi_{1}(x,\xi)\Omega_{1} \nonumber\\
    &+ \int_{\xi}^{x}\Big(K_{1}(\tau,\xi)\Theta_{1}
    +K_{2}(\tau,\xi)\Theta_{2}\Big)\lambda_{1}^{-1}\Phi_{1}(x,\tau)\Omega_{1}d\tau,\label{eq:intpde1_kernels}\\
    &\lambda_{1}\partial_{x} K_{2} - \lambda_{2}\partial_{\xi}K_{2} =-\sigma_{12}K_{1} +  \lambda_{1}^{-1}\Theta_{1}\Phi_{1}(x,\xi)\Omega_{2}\nonumber\\
    &+ \int_{\xi}^{x}\Big( K_{1}(\tau,\xi)\Theta_{1}+K_{2}(\tau,\xi)\Theta_{2}d\tau\Big)\lambda_{2}^{-1}\Phi_{1}(x,\tau)\Omega_{2},\label{eq:intpde2_kernels}\\
    &\lambda_{2}\partial_{x} Q_{1} - \lambda_{1}\partial_{\xi}Q_{1} =\sigma_{21}Q_{2} + \lambda_{2}^{-1}\Theta_{2}\Phi_{2}(x,\xi)\Omega_{1}\nonumber\\
    &+\int_{\xi}^{x}\Big(Q_{1}(\tau,\xi)\Theta_{1} +Q_{2}(\tau,\xi)\Theta_{2}\Big)\lambda_{2}^{-1}\Phi_{2}(\tau,x)\Omega_{1}d\tau,\label{eq:intpde4_kernels}\\
    &\lambda_{2}\partial_{x} Q_{2} + \lambda_{2}\partial_{\xi}Q_{2} = \sigma_{12}Q_{1} + \lambda_{2}^{-1}\Theta_{2}\Phi_{2}(x,\xi)\Omega_{2}\nonumber\\
    &+\int_{\xi}^{x}\Big(Q_{1}(\tau,\xi)\Theta_{1} +Q_{2}(\tau,\xi)\Theta_{2}\Big)\Phi_{2}(\tau,x)\lambda_{2}^{-1}\Omega_{2}d\tau.\label{eq:intpde5_kernels}
\end{align}

These equations, together with boundary conditions \eqref{eq:bc1_kernels}-\eqref{eq:bc2_kernels}, are well-posed as shown in the next Lemma: 
\begin{lemma}
    The PDE system \eqref{eq:intpde1_kernels}-\eqref{eq:intpde5_kernels} with boundary conditions \eqref{eq:bc1_kernels}-\eqref{eq:bc2_kernels} has a unique $\mathcal{C}^{1}(\mathcal{T})$ solution.
\end{lemma}

\begin{proof}
    The well-posedness follows the steps of the proof in the Appendix of \cite{vazquez2011}, with some slight modification to account for the integral terms as in \cite{krstic2008}. \end{proof}
    
    From the previous Lemma and the theory of Volterra integral equations, it follows that the inverse of transformation \eqref{eq_transf1}-\eqref{eq_transf2} always exists and can be defined as
\begin{align}
u=& \alpha  +\int_{x}^{1}L_{1}(x,\xi)\alpha(t,\xi)d\xi  + \int_{x}^{1}L_{2}(x,\xi)\beta(t,\xi)d\xi  \nonumber\\
&+\int_{x}^{1}S(x,\xi)v(t,\xi)d\xi, \label{eq_inv_transf1}\\
p=& \beta+ \int_{x}^{1}M_{1}(x,\xi)\alpha(t,\xi)d\xi +\int_{x}^{1}M_{2}(x,\xi)\beta(t,\xi)d\xi \nonumber\\
&+\int_{x}^{1}E(x,\xi)v(t,\xi)d\xi,\label{eq_inv_transf2}
\end{align}
 where $L_{i}$, $M_{i}$, with $i \in\{ 1, 2\}$, $S = (S_{1}\;\;\; \hdots \;\;\; S_{n})$ and $E = (E_{1}\;\;\; \hdots \;\;\; E_{n})$ are the inverse kernels, which verify  equations similar to \eqref{eq:pde1_kernels}--\eqref{eq:bc3_kernels}.

Using the above results, we state the following Lemma:
\begin{lemma}
Let $U$ be given by the following control law
\begin{multline}
U(t) = - q p(t,0) + \int_{0}^{1}K_{1}(0,\xi)u(t,\xi)d\xi \\
+\int_{0}^{1}K_{2}(0,\xi)p(t,\xi)d\xi + \int_{0}^{1}G(0,\xi)v(t,\xi)d\xi. \label{eq:control_law}
\end{multline}

Then the transformation \eqref{eq_transf1}--\eqref{eq_transf2} maps \eqref{eq_u}--\eqref{bc_u} into the following target system:
\begin{align}
&\partial_{t}\alpha (t,x)=-\lambda_{1}\partial_{x}\alpha (t,x) ,\label{eq:pde1_target}\\
&\partial_{t}\beta (t,x) = \lambda_{2}\partial_{x}\beta (t,x) ,\label{eq:pde2_target}\\
&\partial_{t}v(t,x) = \Omega_{1}\alpha(t,x)+\Omega_{2}\beta(t,x)+\Psi v(t,x)\nonumber\\
&\quad +\int_{x}^{1}N_{1}(x,\xi)\alpha(t,\xi)d\xi +\int_{x}^{1}N_{2}(x,\xi)\beta(t,\xi)d\xi  \nonumber\\
&\quad +\int_{x}^{1}N_{3}(x,\xi)v(t,\xi)d\xi,\label{eq:pde3_target}\\
&\beta (t,1) = \rho \alpha (t,1),\label{eq:bc1_target}\\
&\alpha(t,0) = 0,\label{eq:bc2_target}
\end{align}
\noindent with
\begin{eqnarray*}
N_{j}(x,\xi)&=& \Omega_{1}L_{j}(x,\xi)+\Omega_{2}M_{j}(x,\xi), \qquad \mbox{for } j\in\{1, 2\}, \\
N_{3}(x,\xi)&=& \Omega_{1}S(x,\xi)+\Omega_{2}E(x,\xi),
\end{eqnarray*}
\end{lemma}
\begin{proof}
Differentiating \eqref{eq_transf1}-\eqref{eq_transf2} with respect to time and space, integrating by parts, substituting the resultant expressions into \eqref{eq_u}-\eqref{eq_vn}, and applying \eqref{eq_inv_transf1}-\eqref{eq_inv_transf2}, we obtain \eqref{eq:pde1_target}--\eqref{eq:pde3_target}.

Evaluating \eqref{eq_transf1} at $x=0$, substituting it into \eqref{bc_u} and using \eqref{eq:control_law} we obtain \eqref{eq:bc1_target}. Finally, evaluating \eqref{eq_transf2} for $x=1$, substituting it into \eqref{bc_p} and using \eqref{eq_inv_transf2} we obtain \eqref{eq:bc1_target}.
\end{proof}

\subsection{Stability of the target system}

The stability properties of the target system \eqref{eq:pde1_target}-\eqref{eq:bc2_target} are proved in the following Lemma:

\begin{lemma}
    The zero equilibrium of system \eqref{eq:pde1_target}-\eqref{eq:bc2_target} is exponentially stable in the $L_{2}$ sense.
\end{lemma}
\begin{proof}
    Consider the following Lyapunov functional
    \begin{align}
        V(t) =& \int_{0}^{1}\left(\frac{A}{\lambda_{1}}\mathrm{e}^{-\mu x} \alpha^{2}(t,x) +\frac{B}{\lambda_{2}}\mathrm{e}^{\mu x}\beta^{2}(t,x)\right) dx \nonumber \\ 
        &\displaystyle +\frac{1}{2}\int_{0}^{1}v^{T}(t,x)P(x)v(t,x)dx,
    \end{align}
    \noindent where $P(x) = \mathrm{e}^{-\vartheta x}I_{n\times n}$,in which $I_{n\times n}$ stands for the $n\times n$ identity matrix, and $A$, $B$, $\mu$, and $\vartheta$ are constants to be defined. Differentiating $V$ with respect to time yields
    \begin{align}
        \dot{V}(t) =& 2\int_{0}^{1}\bigg(-A\mathrm{e}^{-\mu x} \alpha(t,x)\partial_{x}\alpha(t,x) \nonumber\\
        &+B\mathrm{e}^{\mu x}\beta(t,x)\partial_{x}\beta(t,x)\bigg) dx \nonumber\\
        &+ \int_{0}^{1}v^{T}(t,x)P(x)\bigg(\Omega_{1}\alpha(t,x) \nonumber\\
        &+\Omega_{2}\beta(t,x) +\Psi v(t,x)\bigg)dx \nonumber \\
        &+ \int_{0}^{1}v^{T}(t,x)P(x)\left(\int_{x}^{1}\bigg(N_{1}(x,s)\alpha(t,s)\right.\nonumber\\
        &\left.+N_{2}(x,s)\beta(t,s) + N_{3}(x,s)v(t,s)\bigg)ds \right)dx.\label{eq:derivative_V}
    \end{align}

    Integrating by parts the first two terms in the right-hand side of \eqref{eq:derivative_V} and plugging the boundary conditions \eqref{eq:bc1_target}-\eqref{eq:bc2_target} yields
    \begin{align}
        & \int_{0}^{1}\left(-A\mathrm{e}^{-\mu x} \alpha(t,x)\partial_{x}\alpha(t,x)+B\mathrm{e}^{\mu x}\beta(t,x)\partial_{x}\beta(t,x)\right) dx\nonumber\\
        &= -\mu\int_{0}^{1}\bigg(A\mathrm{e}^{-\mu x}\alpha^{2}(t,x)+B\mathrm{e}^{\mu x}\beta^{2}(t,x)\bigg)dx  \nonumber\\
        &\quad+ (B\rho^{2}\mathrm{e}^{\mu}-A\mathrm{e}^{-\mu})\alpha^{2}(t,1) -B\beta^{2}(t,0). \label{eq:bound1}
    \end{align}

 To make further progress, we will now compute an upper bound for the rest of the terms in the right-hand side of \eqref{eq:derivative_V}.
 
 By using the Young's inequality and considering the fact that $\mathrm{e}^{-\vartheta x}\leq 1\leq \mathrm{e}^{\mu }\mathrm{e}^{-\mu x}$ and $\mathrm{e}^{-\delta x}\leq1 \leq \mathrm{e}^{\mu x}$, we have that
    \begin{align}
    &\int_{0}^{1}v^{T}(t,x)P(x)\bigg(\Omega_{1}\alpha(t,x)+\Omega_{2}\beta(t,x)\bigg)dx\nonumber\\
    \leq &\frac{2\mathrm{e}^{\mu }}{\rho(\Psi)}\Omega_{1}^{T}\Omega_{1}\int_{0}^{1}\mathrm{e}^{-\mu x}\alpha^{2}(t,x)dx \\
    &+\frac{2}{\rho(\Psi)}\Omega_{2}^{T}\Omega_{2}\int_{0}^{1}\mathrm{e}^{\mu x}\beta^{2}(t,x)dx\nonumber \\
    &+ \frac{\rho(\Psi)}{4}\int_{0}^{1}v^{T}(t,x)P(x) v(t,x)dx,\label{eq:bound2}
    \end{align}
    \noindent where $\rho(\Phi)$ is the spectral radius of $\Psi$.

Now, define $\overline{N}_{i} = \|N_{i}(x,y)\|_{\infty}$, for $i\in\{1,2,3\}$. Then, using the Cauchy-Schwarz and Young's inequalities, we get
    \begin{align}
        & \int_{0}^{1}\int_{x}^{1}v^{T}(t,x)P(x) \left(N_{1}(x,s)\alpha(t,s) +N_{2}(x,s)\beta(t,s) \right)ds dx\nonumber\\
        &\leq  \frac{2n\overline{N}_1^{2}\mathrm{e}^{\mu}}{\rho(\Psi)} \int_{0}^{1}\mathrm{e}^{-\mu x}\alpha^{2}(t,x)dx + \frac{2n \overline{N}_2^{2}}{\rho(\Psi)} \int_{0}^{1}\mathrm{e}^{\mu x}\beta^{2}(t,x)dx \nonumber \\
        &\quad+ \frac{\rho(\Psi)}{4}\int_{0}^{1}v^{T}(t,x)P(x)v(t,x)dx.  \label{eq:bound_V2}
    \end{align}

Finally, 
\begin{align*}
    &\int_{1}^{x} |N_{3}(x,s) v (t,s)| ds \\
    &\leq\overline{N}_3\mathrm{e}^{\vartheta/2}  \frac{\mathrm{e}^{\vartheta x/2}}{\sqrt{\vartheta}}  \sqrt{ \int_{0}^{1}v^{T}(t,s)P(s)v(t,s)ds},
\end{align*}
where at the last step, Cauchy-Schwarz and Young's inequalities were again applied.

Therefore we reach
\begin{align*}
    \dot{V}\leq& -\left(A\sigma-\frac{2(\Omega_{1}^{T}\Omega_{1}+\overline{N}_1^2) \mathrm{e}^{\sigma}}{\rho(\Phi)}\right) \int_{0}^{1}\mathrm{e}^{-\mu x}\alpha^{2}(t,x)dx   \nonumber \\
    &-\left(B\sigma- \frac{2(\Omega_{2}^{T}\Omega_{2}+\overline{N}_2^2 )}{\rho(\Phi)}\right) \int_{0}^{1}\mathrm{e}^{\mu x}\beta^{2}(t,x)dx \nonumber \\
    & -\left( \frac{\rho(\Phi)}{2}-\frac{ \overline{N}_3}{\sqrt{\vartheta}} \right)  \int_{0}^{1}v^{T}(t,x)P(x)v(t,x)dx\nonumber \\ 
    & + (Aq^{2}-B)\beta^{2}(t,0).
\end{align*}

Choosing 
\begin{align*}
     A=&\mathrm{e^\sigma}, \quad B=Aq^2+1, \quad     \vartheta=\frac{ 16 \overline{N}_3^2}{\rho(\Phi)^2},\\
    \mu=&\max\left\{
    \frac{2 (\Omega_{1}^{T}\Omega_{1}+\overline{N}_1^2 )}{\rho(\Phi)},
    \frac{2 (\Omega_{2}^{T}\Omega_{2}+\overline{N}_2^2 )}{\rho(\Phi)}
    \right\}+1,
\end{align*}
we get
\begin{align*}
    \dot{V}\leq & - \int_{0}^{1}\mathrm{e}^{-\mu x}\alpha^{2}(t,x)dx  -\int_{0}^{1}\mathrm{e}^{\mu x}\beta^{2}(t,x)dx \\
    &-\frac{\rho(\Phi)}{4}  \int_{0}^{1}v^{T}(t,x)P(x)v(t,x)(t,x)dx \\
    \leq& -K V
\end{align*}
for $K=\min\left\{\frac{2\lambda_1}{A},\frac{2\lambda_2}{B},\frac{\rho(\Phi)}{2} \right\}>0$, thus proving exponential stability of the equilibrium $\alpha\equiv\beta\equiv v\equiv0$.
\end{proof} 

\section{Numerical simulations}\label{section:simulations}

In this section, we present numerical simulations of system \eqref{eq_u}-\eqref{bc_p} with the proposed control law \eqref{eq:control_law} considering $n=2$. The parameters were chosen as $q=-0.7$, $\rho=0.5$, $\sigma_{12}= 2.5$, $\sigma_{21}=-3.5$, $\theta_{11}=0.25$, $\theta_{12}=0.1$, $\theta_{21}=0.25$, $\theta_{22}=-0.1$, $\omega_{11}=0.3$, $\omega_{12}=0.8$, $\omega_{21}=-0.65$, $\omega_{22}=0.3$, $\psi_{11}=-1.5$, $\psi_{12}=2$, $\psi_{21}=-1$, $\psi_{22}=-2$, $\lambda_{1}= 1.25$ and $\lambda_{2}=0.9$, which corresponds to an open-loop unstable system. The finite differences method was employed in MATLAB to compute the states of the system and solve the kernel PDEs.

\begin{figure*}[th!]
    \centering
    \includegraphics[width=60mm]{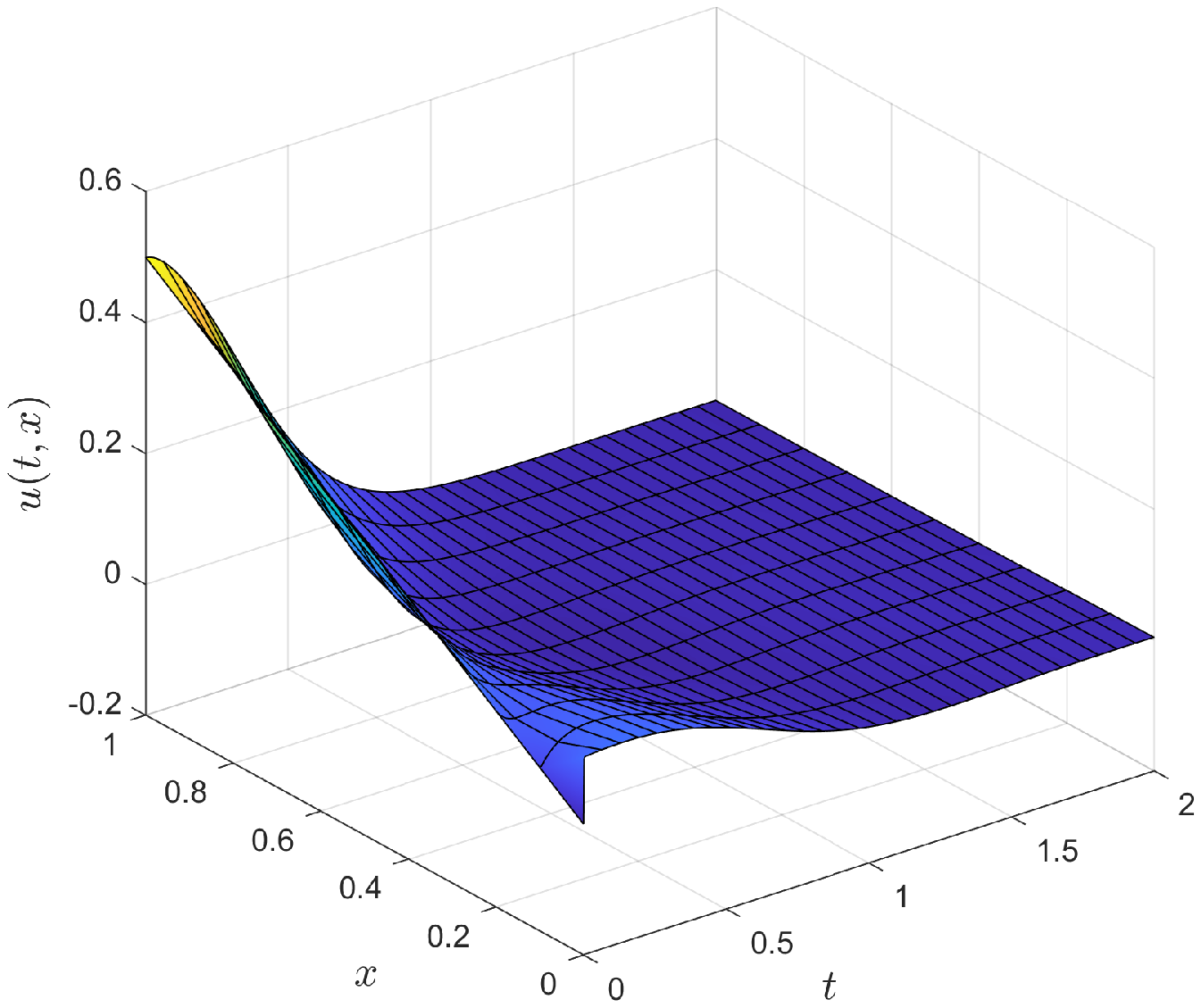}
    \includegraphics[width=60mm]{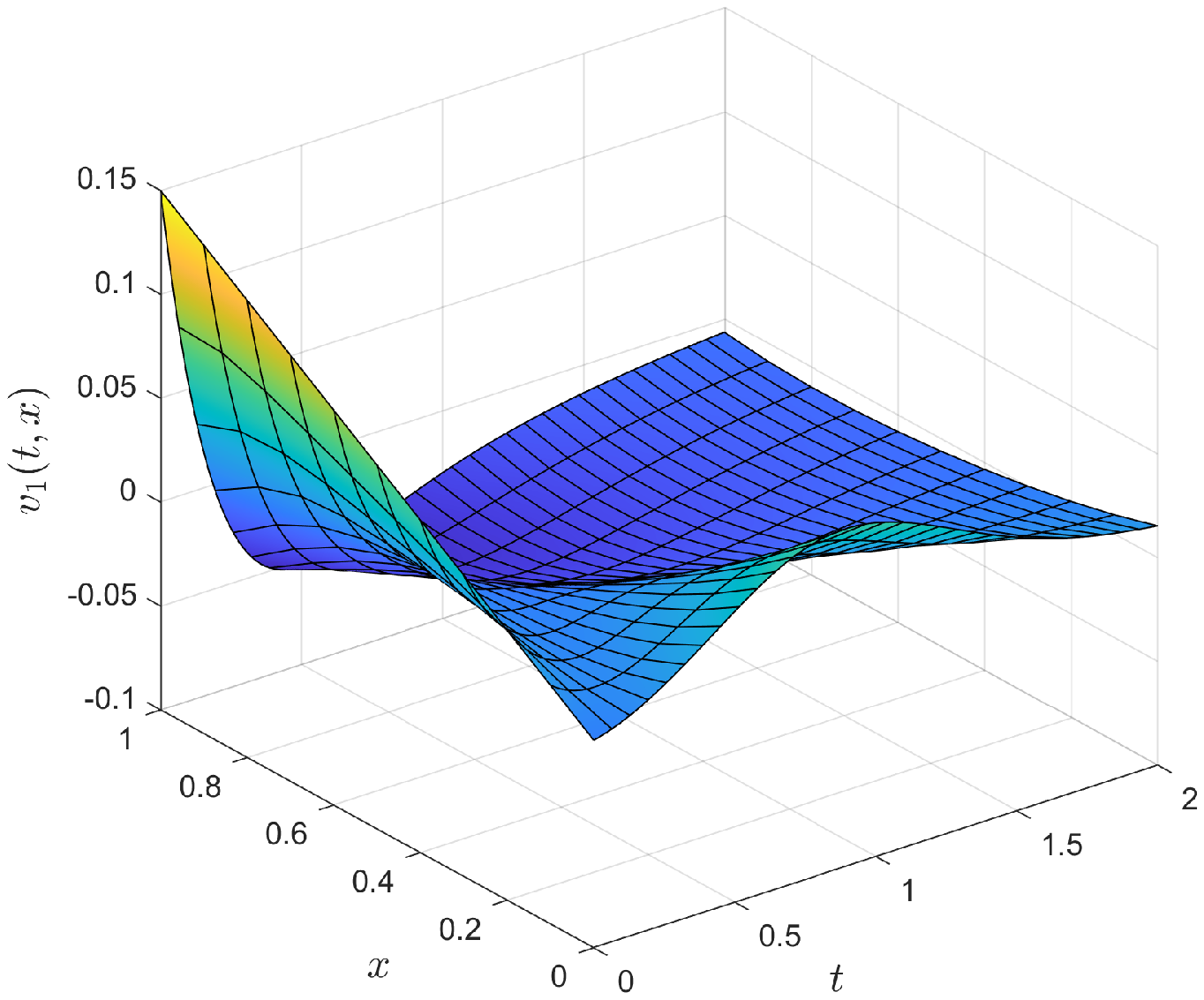}\\
    \includegraphics[width=60mm]{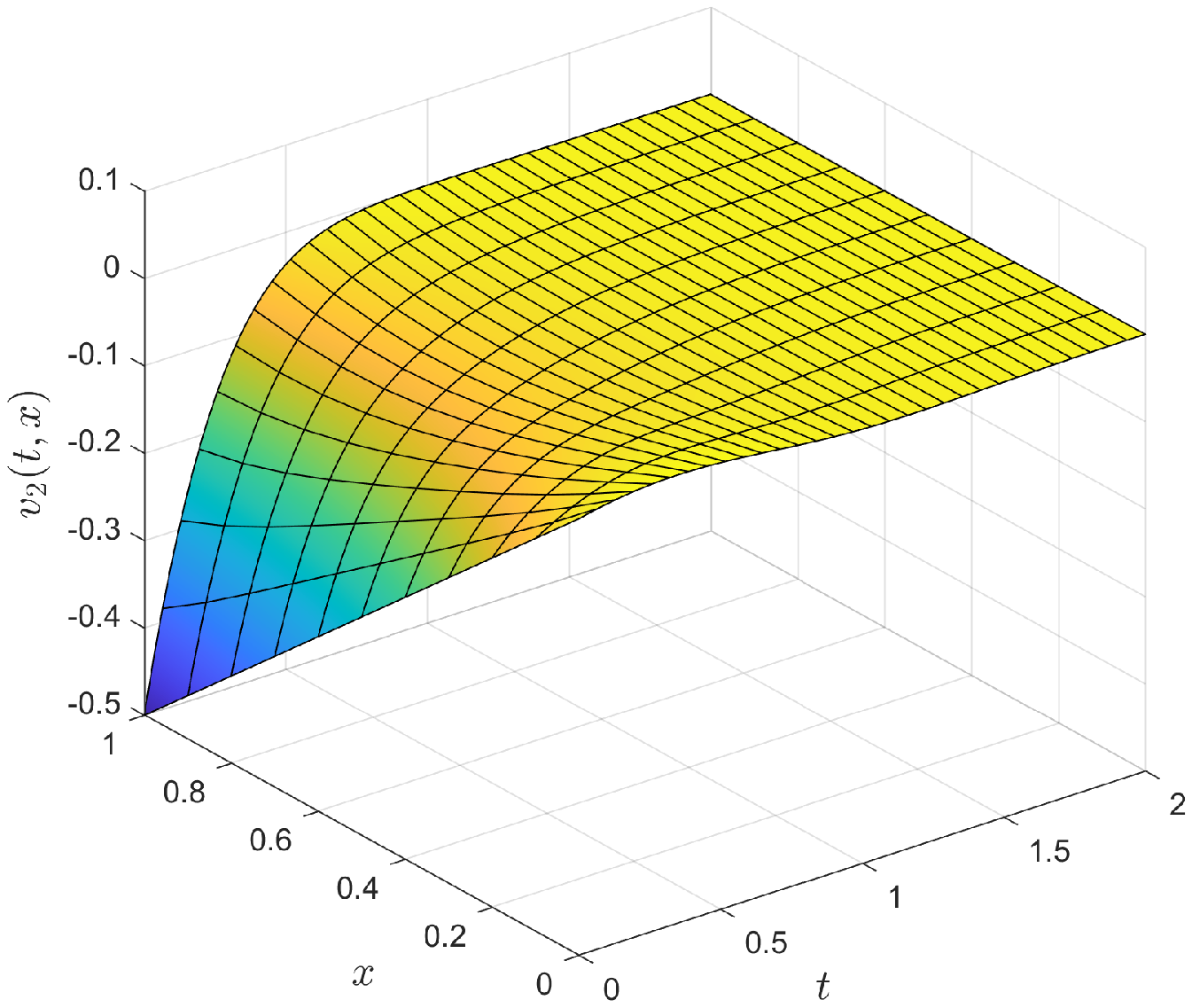}
    \includegraphics[width=60mm]{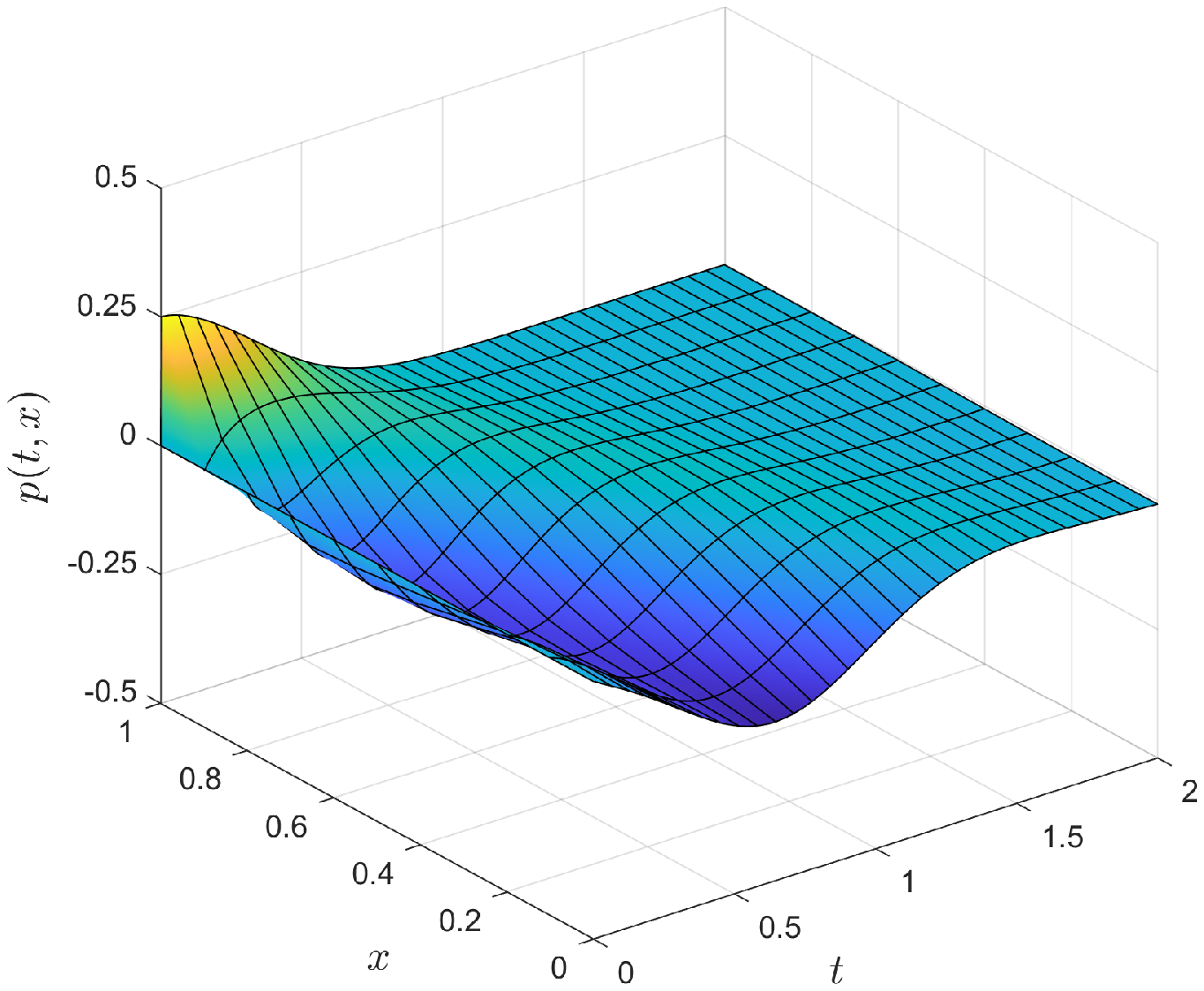}
    \caption{Distributed states evolution as a function of time and space.}
    \label{fig:closed_loop_states}
\end{figure*}

\begin{figure}[h!]
    \centering
    \includegraphics[width=60mm]{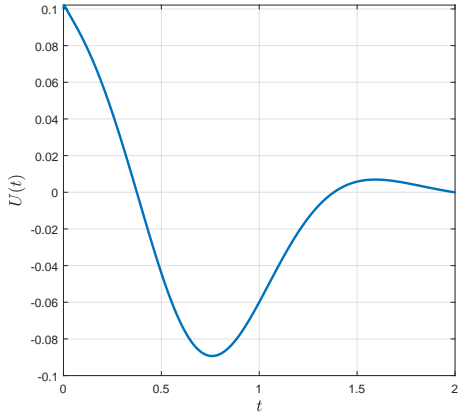}
    \caption{State feedback control law as a function of time.}
    \label{fig:control_law}
\end{figure}

Figures \ref{fig:closed_loop_states} and \ref{fig:control_law}, show the closed-loop states and the control signal, respectively. As can be seen in Figures  \ref{fig:closed_loop_states} and \ref{fig:control_law}, the system states decay to zero after an initial transient.

\section{Conclusion}\label{section:conclusions}

In this work, we introduced the state feedback stabilization of a class of hyperbolic systems containing an atachic subsystem, with zero characteristic velocities, which we denote as $(1+n+1)\times(1+n+1)$ systems. We showed that stabilizability requires that the atachic subsystem  be asymptotically stable. Under this condition, we proved that the boundary feedback stabilization of these systems is possible by applying the backstepping methodology to guarantee the closed-loop exponential stability in the $\mathcal{L}^{2}$ sense. Interestingly, we employ the invertible Volterra transformation only for the PDEs with nonzero characteristic speeds, leaving the atachic subsystem  unaltered in the target system but making it ISS with respect to the counterconvecting nonzero-speed states. Although only the state feedback controller was studied in this article, a Luenberguer-type state observer with boundary measurements of states with non-zero characteristic speeds can be designed to obtain the associated output feedback controller. Future discussions should include general heterodirectional hyperbolic PDEs or networks of systems of hyperbolic balance laws coupled with ODEs. 

\appendix

\subsection{Proof of Proposition 2.1}\label{app:proof_prop1}

Let arbitrary $(u_{0},v_{0})\in \mathcal{H}^{1}(0,1)\times \mathcal{H}^{1}(0,1)\setminus S$ be given. For all $x\in[0,1]$, define the functions
\begin{align}
    \varphi (x) &\triangleq \psi\int_{0}^{x}v_{0}(s)ds + \omega\int_{0}^{x}u_{0}(s)ds +  \lambda v_{0}(x),\label{eq:prop_func_phi}\\
    h(x) &\triangleq \varphi (x) - K,\label{eq:prop_func_h}
\end{align}
\noindent where
\begin{align}
    K \triangleq \int_{0}^{1}\varphi (x)dx.\label{eq:prop_func_K}
\end{align}

Definitions \eqref{eq:linear_subspace}, \eqref{eq:prop_func_phi}-\eqref{eq:prop_func_K} give the implication 
\begin{align}
    h(x)\equiv 0 \Rightarrow (u_{0},v_{0})\in S.\label{eq:prop_implication_h}
\end{align}

Therefore, since $(u_{0},v_{0})\notin S$ it follows from \eqref{eq:prop_implication_h} that $h\in \mathcal{H}^{1}(0,1)$ is a non-identically zero function, i.e., $h\neq 0$. Moreover, definitions \eqref{eq:prop_func_h}-\eqref{eq:prop_func_K} imply that $\int_{0}^{1}h (x)dx =0$. Since the Cauchy-Schwarz inequality in $\mathcal{L}^{2}(0,1)$ holds as an equality if and only if the functions are linearly dependent in $\mathcal{L}^{2}(0,1)$, we obtain from  \eqref{eq:linear_subspace}, \eqref{eq:prop_func_phi}-\eqref{eq:prop_func_K} that
\begin{align}
    \|\varphi\|_{2} = |K| \Leftrightarrow (u_{0},v_{0})\in S.\label{eq:prop_norm2_phi}
\end{align}

Consequently, it follows from \eqref{eq:prop_norm2_phi} that
\begin{align}
    &\int_{0}^{1}h(x)\left(\psi\int_{0}^{x}v_{0}ds+ \omega\int_{0}^{x}u_{0}(s)ds\right. \nonumber\\
    & \qquad+  \lambda v_{0}(x) \bigg)dx = \|\varphi\|_{2}^{2} - K^{2} > 0.\label{eq:prop_inequality}
\end{align}

We show next that for every input $U\in \mathcal{C}^{1}(\mathbb{R}_{+})$ with $U(0)=u_{0}(0)$ the unique solution $(u,v)\in \left(\mathcal{C}^{0}(\mathbb{R}_{+};\mathcal{H}^{1}(0,1))\cap\right.$ $\left. \mathcal{C}^{1}(\mathbb{R}_{+};\mathcal{L}^{2}(0,1))\right)^{2}$ of the initial-boundary value problem \eqref{eq:gen_simplified_u}-\eqref{eq:gen_ini_simplified} does not satisfy $\lim_{t\rightarrow\infty}(\|u[t]\|_{2}) = \lim_{t\rightarrow\infty}(\|v[t]\|_{2}) = 0.$ 

By contradiction, suppose that there exists an input $U\in \mathcal{C}^{1}(\mathbb{R}_{+})$ with $U(0)=u_{0}(0)$ for which the unique solution $(u,v)\in \left(\mathcal{C}^{0}(\mathbb{R}_{+};\mathcal{H}^{1}(0,1))\cap\right.$ $\left. \mathcal{C}^{1}(\mathbb{R}_{+};\mathcal{L}^{2}(0,1))\right)^{2}$ of the initial-boundary value problem \eqref{eq:gen_simplified_u}-\eqref{eq:gen_ini_simplified} satisfies $\lim_{t\rightarrow\infty}(\|u[t]\|_{2}) = \lim_{t\rightarrow\infty}(\|v[t]\|_{2}) = 0$. Define the functional
\begin{align}
    R(t) \triangleq& \int_{0}^{1}h(x)\left(\psi\int_{0}^{x}v(t,s)ds \right.\nonumber\\
    &+ \left.\omega\int_{0}^{x}u(t,s) ds +  \lambda v(t,x)\right)dx,\label{eq:prop_functional_R}
\end{align}
\noindent for $t\geq 0$.

Since $(u,v)\in \left( \mathcal{C}^{1}(\mathbb{R}_{+};\mathcal{L}^{2}(0,1))\right)^{2}$, it follows that $R\in \mathcal{C}^{1}(\mathbb{R}_{+})$. Using \eqref{eq:bc_simplified}, \eqref{eq:prop_functional_R}, and the fact that  $\int_{0}^{1}h(x)dx=0$, we get, after substituting \eqref{eq:gen_simplified_u}-\eqref{eq:gen_bc_simplified} and integrating by parts, for all $t\geq 0$:
\begin{align*}
    \dot{R}(t) &= \int_{0}^{1}h(x)\left(\psi^{2}\int_{0}^{x}v(t,s)ds + \psi\omega\right.\\
    &\times\left. \int_{0}^{x}u(t,s)ds +  \lambda\psi v(t,x) +  \lambda\omega u(t,0) \right)dx\\
    &= \psi R(t) +  \lambda\omega U(t)\int_{0}^{1}h(x)dx = \psi R(t).
\end{align*}

From the above expression, it follows that
\begin{align}
    R^{2}(t)= \mathrm{e}^{2\psi t}R^{2}(0),\label{eq:solution_R}
\end{align}
\noindent for all $t\geq 0$.

Notice that definitions \eqref{eq:prop_inequality} and \eqref{eq:prop_functional_R} imply that $R(0)\neq 0$ and since $\psi \geq 0$, we get from \eqref{eq:solution_R} that
\begin{align}
    R^{2}(t)\geq R^{2}(0)>0,\label{eq:inequality_R}
\end{align}
\noindent for all $t\geq 0$. Using definition \eqref{eq:prop_functional_R} and the Cauchy-Schwarz and Holder inequalities, we obtain the estimate:
\begin{align}
    R^{2}(t)\leq& 2(\psi + \lambda)^{2}\|h\|_{\infty}^{2}\|v[t]\|_{2}^{2}
                &+ 2|\omega |\|h\|_{\infty}^{2}\|u[t]\|_{2}^{2},\label{eq:estimate_R}
\end{align}
\noindent for all $t\geq 0$. However, inequalities \eqref{eq:inequality_R}-\eqref{eq:estimate_R} contradict the fact that $\lim_{t\rightarrow\infty}(\|u[t]\|_{2}) = \lim_{t\rightarrow\infty}(\|v[t]\|_{2}) = 0$. 

Therefore, for every input $U\in \mathcal{C}^{1}(\mathbb{R}_{+})$ with $U(0)=u_{0}(0)$ the corresponding unique solution $(u,v)\in \left(\mathcal{C}^{0}(\mathbb{R}_{+};\mathcal{H}^{1}(0,1))\cap\right.$ $\left. \mathcal{C}^{1}(\mathbb{R}_{+};\mathcal{L}^{2}(0,1))\right)^{2}$ of the initial-boundary value problem \eqref{eq:gen_simplified_u}-\eqref{eq:gen_ini_simplified} does not satisfy $\lim_{t\rightarrow\infty}(\|u\|_{2}) = \lim_{t\rightarrow\infty}(\|v\|_{2}) = 0.$  

\subsection{Proof of Proposition 2.2}\label{app:proof_prop2}

Define
\begin{align}
    w(t,x) =& v(t,x)- \mathrm{e}^{- \lambda^{-1}\psi x}v(t,0)\nonumber\\
&+  \lambda^{-1} \omega \int_{0}^{x}\mathrm{e}^{- \lambda^{-1}\psi (x-s)}u(t,s)ds,\label{eq:definition_w}
\end{align}
\noindent for all $t\geq 0$ and $x\in[0,1]$.

Note that \eqref{eq:gen_simplified_u}-\eqref{eq:gen_bc_simplified} imply the following equation:
\begin{align}
    \partial_{t}w = \psi w.\label{eq:pde_w}
\end{align}

Definitions \eqref{eq:linear_subspace} and \eqref{eq:definition_w} give the following equivalance for all $t\geq 0$:
\begin{align}
    (u[t],v[t])\in S \Leftrightarrow w[t]=0,
\end{align}
\noindent while \eqref{eq:pde_w} gives the implication $w[0]=0\rightarrow w[t]=0$ for all $t\geq 0$. The proof is complete.

\subsection{Proof of Proposition 2.3}\label{app:proof_prop3}

Let arbitrary $(u_{0},v_{0})\in \mathcal{H}^{1}(0,1)\times \mathcal{H}^{1}(0,1)\setminus S$ be given. We show next that for every input $U\in \mathcal{C}^{1}(\mathbb{R}_{+})$ with $U(0)=u_{0}(0)$ the unique solution $(u,v)\in \left(\mathcal{C}^{0}(\mathbb{R}_{+};\mathcal{H}^{1}(0,1))\cap\right.$ $\left. \mathcal{C}^{1}(\mathbb{R}_{+};\mathcal{L}^{2}(0,1))\right)^{2}$ of the initial-boundary value problem \eqref{eq:gen_simplified_u}-\eqref{eq:gen_ini_simplified} does not satisfy $\lim_{t\rightarrow\infty}v(t,x) = 0$ for all $x\in[0,1]$.

For the solution of \eqref{eq:gen_simplified_u}-\eqref{eq:gen_ini_simplified} the following formulas are valid for $t\geq 0$ and $x\in[0,1]$:
\begin{align}
    u(t,x) &= \left\{ \begin{array}{ll}
         u_{0}(x- \lambda t), & 0\leq t \leq  \lambda^{-1} x, \\
         U(t- \lambda^{-1}x),& t>  \lambda^{-1} x, 
    \end{array}\right. \label{eq:solution_u}\\
    v(t,x)&=\mathrm{e}^{\psi t}v_{0}(x) + \omega \int_{0}^{x}\mathrm{e}^{\psi (t-s)}u(s,x)ds.\label{eq:solution_v}
\end{align}

It follows from \eqref{eq:gen_simplified_u}-\eqref{eq:gen_bc_simplified} that the following equation holds for all $x\in[0,1]$ and $t> \lambda^{-1}x$:
\begin{multline}
    \frac{d}{dt}\left(v(t- \lambda^{-1}x,0) - v(t,x) \right) = \\ \psi \left(v(t- \lambda^{-1}x,0)-v(t,x). \right)\label{eq:difference_v}
\end{multline}

Using \eqref{eq:difference_v} and continuity of $v$, we get for all $x\in[0,1]$ and $t\geq  \lambda^{-1}x$:
\begin{multline}
    v(t- \lambda^{-1}x,0) - v(t,x) =\\ \mathrm{e}^{\psi (t- \lambda^{-1}x)}\left(v_{0}(0) - v( \lambda^{-1}x,x) \right)\label{eq:difference2_v}
\end{multline}

Since $\psi \geq 0$, we get from \eqref{eq:difference2_v} for all $x\in[0,1]$ and $t\geq  \lambda^{-1}x$:
\begin{multline}
    \left|v(t- \lambda^{-1}x,0) - v(t,x)\right| \geq\\
    \mathrm{e}^{\psi (t- \lambda^{-1}x)} \left| v_{0}(0) - v( \lambda^{-1}x,x) \right|\label{eq:inequality_difference_v}
\end{multline}

Inequality \eqref{eq:inequality_difference_v} combined with the fact that $\lim_{t\rightarrow \infty}v(t,x)=0$ for all $x\in[0,1]$ implies that $v( \lambda^{-1}x,x) = v_{0}(0) $ for all $x\in[0,1]$. This equation combined with \eqref{eq:solution_v} gives for all $x\in[0,1]$:
\begin{align}
   v_{0}(0)  =& \mathrm{e}^{\psi  \lambda^{-1}x}v_{0}(x) \nonumber\\
   &+ \omega  \int_{0}^{ \lambda^{-1}x}\mathrm{e}^{\psi ( \lambda^{-1}x-s)} u(s,x)ds.\label{eq:prop3}
\end{align}

Using \eqref{eq:solution_u} and \eqref{eq:prop3} we get for all $x\in[0,1]$:
\begin{align}
    v_{0}(0)  = v_{0}(x) + \omega \int_{0}^{ \lambda^{-1}x}\mathrm{e}^{-\psi s}u_{0}(x- \lambda s)ds.\label{eq:prop3_aux2}
\end{align}

However, \eqref{eq:prop3_aux2} in conjunction with definitions \eqref{eq:linear_subspace}-\eqref{eq:operator_P} imply that $(u_{0},v_{0})\in S$; a contradiction.

Therefore, for every input $U\in \mathcal{C}^{1}(\mathbb{R}_{+})$ with $U(0)=u_{0}(0)$ the corresponding unique solution $(u,v)\in \left(\mathcal{C}^{0}(\mathbb{R}_{+};\mathcal{H}^{1}(0,1))\cap\right.$ $\left. \mathcal{C}^{1}(\mathbb{R}_{+};\mathcal{L}^{2}(0,1))\right)^{2}$ of the initial-boundary value problem \eqref{eq:gen_simplified_u}-\eqref{eq:gen_ini_simplified} does not satisfy $\lim_{t\rightarrow\infty}v(t,x) = 0.$ for all $x\in[0,1]$. 

\bibliographystyle{IEEEtran}
\bibliography{bibliography}

%








\end{document}